\documentclass[smallextended]{amsart}
\usepackage{amssymb}
 \usepackage{graphicx}
\usepackage{amscd}
\usepackage{amsmath}
\usepackage{amsmath}
\usepackage{enumitem}      

\usepackage{tikz-cd}

\usepackage{amsfonts,latexsym,amstext}
\usepackage [T1]{fontenc}
\usepackage [latin1]{inputenc}
\usepackage{color}
\usepackage[all]{xy}
\usepackage{cancel}
\usepackage{enumitem} 

\usepackage{tikz-cd}
\usetikzlibrary{arrows}

\newcommand{\Pol}[2]{\ensuremath{{\mathcal P}(^{#1}#2)}}
\newcommand{\Poll}[3]{\ensuremath{{\mathcal P}(^{#1}#2, #3)}}

\newcommand{\tensor}[2]{\ensuremath{{#1}\hat{\otimes}_{\pi}\cdots
\hat{\otimes}_{\pi}{#2}}}

\newcommand{\cartesian}[3]{\ensuremath{#1_{#2} \times \cdots \times #1_{#3}}} 
\newcommand{\Tensor}[4]{\ensuremath{#1_{#2}}{\otimes}_{#4}\cdots{\otimes}_{#4}#1_{#3}} 
\newcommand{\CTensor}[4]{\ensuremath{#1_{#2}}{\hat{\otimes}}_{#4}\cdots{\hat{\otimes}}_{#4}#1_{#3}}

\newcommand{\N}{{\mathbb N}}
\newcommand{\R}{{\mathbb R}}
\newcommand{\C}{{\mathbb C}}
\newcommand{\funcion}[3]{\ensuremath{#1\colon#2 \longrightarrow #3}}

\newcommand{\norma}[1]{\ensuremath{\left\|#1\right\|}}
\newcommand{\normasub}[2]{\ensuremath{\|#1{\|}_{#2}}}

\newcommand{\spaces}[3]{\ensuremath{#1_{#2}, \ldots , #1_{#3}}} 

\newtheorem{theorem}{Theorem}[section]
\newtheorem{definition}[theorem]{Definition}
\newtheorem{question}[theorem]{Question}
\newtheorem{lemma}[theorem]{Lemma}
\newtheorem{proposition}[theorem]{Proposition}
\newtheorem{corollary}[theorem]{Corollary}
\newtheorem{example}[theorem]{Example}
\newtheorem{remark}[theorem]{Remark}

\newenvironment{pf}{\noindent{\it Proof.}}{\\}


\begin{document}

\title{\textbf Lipschitz $p$-summing  multilinear operators}

\author{Jorge C.  Angulo-L\'opez}
\address{Universidad An\'ahuac Mayab; Carretera M\'erida-Progreso km 15.5 AP. 96 Cordemex CP. 97310 M\'erida, Yucat\'an, M\'exico}.

 \email{jorge.angulo@anahuac.mx}

\author{Maite Fern\'{a}ndez-Unzueta}
\address{Centro de Investigaci\'{o}n en Matem\'{a}ticas (Cimat), A.P. 402 Guanajuato, Gto., M\'{e}xico\\
orcid{ 0000-0002-8321-4877}} \email{maite@cimat.mx}

\thanks{J.C. Angulo-L\'opez was supported by a scholarship and  M. Fern\'andez-Unzueta by  project 284110, both  from CONACyT}


%
%
%
%
%
%
%
%
%
%
%
%
%
%
%
%
%

\begin{abstract}
 We apply  the geometric approach provided  by  $\Sigma$-operators to  develop a theory of $p$-summability  for multilinear operators. In this way,  we introduce the notion of Lipschitz $p$-summing multilinear operators and show that it  is  consistent with a general panorama of generalization: Namely, they satisfy Pietsch-type domination and factorization   theorems and  generalizations of the inclusion Theorem,  Grothendieck's coincidence  Theorems,  the  weak Dvoretsky-Rogers Theorem and  a Lindenstrauss-Pe\l czy\'nsky Theorem.  We also characterize  this new class in tensorial terms by means of a
 Chevet-Saphar-type tensor norm. Moreover, we introduce the notion of Dunford-Pettis multilinear operators. With them,  we    characterize  when a projective tensor product  contains $\ell_1$.
   Relations between Lipschitz  $p$-summing multilinear operators with Dunford-Pettis and  Hilbert-Schmidt multilinear
operators are  given.
\end{abstract}

\maketitle

\section{Introduction}\label{sec: intro and pre}

   The relevant role that absolutely $p$-summing operators play in the theory of Banach spaces has motivated the development of analogous notions for  classes of mappings other than linear bounded operators. This is the case of    completely
   $p$-summing operators,  studied by G. Pisier in \cite{Pisier Non-comm p-sum}, and the case of  Lipschitz $p$-summing maps,  studied by J.Farmer and W.B. Johnson in \cite{FJ09}.
This is  also the case
  of bounded multilinear and polynomial  mappings. In this context, however,   a variety  of notions of $p$-summability have  appeared.  In general, they are not equivalent to each other.  Among them, we find   \cite{BomPerVil},  \cite{Dimant03},  \cite{Matos03},   \cite{Pietsch}.    Relations among these notions may be found  in   \cite{CaliPelle}, \cite{PelRib} and \cite{Per-Gar}.
The   current  interest in studying $p$-summability in such a  non-linear context has been stimulated also  by the search of
  Bohnenblust-Hille or   Littlewood-type estimates (e.g. \cite{AlbBayPelSeo}, \cite{BlBoPelRu}, \cite{DefSev}, \cite{DimSev}).

In this paper we  introduce and develop the  notion of   $p$-summability for  multilinear operators that  stems from  regarding, via a natural isomorphism,  multilinear operators  as homogenous  mappings  on a certain   subset of the tensor product, i.e.  as $\Sigma$-operators.  This approach has been successfully applied  to the case of  operators that factor through a Hilbert space \cite{F-UG-H}, as well as to the case of  $(p,q)$-dominated operators \cite{F-UG-H 2}.  Multilinear operators of this type will be called Lipschitz $p$-summing.   Our first result is   that  the triad of equivalent formulations of the linear $p$-summability property holds: the local definition,  the Pietsch's  Domination theorem  and   the Pietsch's    Factorization theorem (see Theorem \ref{thm: equiv p summ multi op}).

 To provide a general framework for this class of multilinear mappings, we prove multilinear generalizations of several other fundamental results of the linear theory.
Specifically, we prove that multilinear versions of the  Inclusion Theorem (Proposition  \ref{prop: inclusion thm}), as well as Grothendieck's coincidence  Theorems  (Corollaries \ref{coro: Groth Thm} and \ref{coro pi1 = pi2}) hold.
    Multilinear results generalizing   the weak Dvoretzky-Rogers'  and Lindenstrauss-Pe\l czy\'nsky's theorems are  also proved (Corollary \ref{coro: teorema debil de Dvoretzky-Rogers} and Corollary \ref{coro l1 Hilbert}, respectively).
    In addition, we describe
    Lipschitz $p$-summing multilinear operators   as elements in the dual of a tensor space, whose norm generalizes the Chevet-Saphar crossnorm  (see Theorems \ref{dual representation} and \ref{thm: equiv p-Lips norm}).

     In general, Lipschitz $p$-summing multilinear operators are not comparable with compact operators. However,
     as it happens in the linear case,  they transform some sequences with a weak-Cauchy behaviour into norm-convergent sequences (see Proposition \ref{prop: psum DP}).  Multilinear operators satisfying this last property will be called Dunford-Pettis operators. They  are studied in Section \ref{sec: DP}.

    In all our considerations  above, we have been implicitly using the projective tensor product. In Section \ref{sec: Hil-Sch} we  briefly see how to develop the theory of Lipschitz $p$-summing operators with respect to any reasonable crossnorm $\beta$.  We apply this to see that  for every $1 \leq p < \infty$, Hilbert-Schmidt multilinear operators coincide with  the class of  {\sl Lipschitz $p$-summing multilinear operators with respect to the Hilbert-Schmidt tensor product},  which generalizes a classical linear result proved in \cite{Pelcz66-67} (see Theorem \ref{thm: Hilb-Schm}).

 Finally,  in Section \ref{sec: comparison} we introduce   Lipschitz $p$-summing homogeneous polynomials and relate these notions with other non-linear generalizations of absolutely $p$-summing operators that have appeared in the literature.

\section{Multilinear mappings studied through $\Sigma$-operators.}\label{sec: sigma-op}

Throughout this paper $X,X_1,\ldots, X_n$ and $Y$ will be  Banach spaces
over the same field $\R$ or $\C$, and  $B_X$ will be   the closed unit ball of a space   $X$.  The projective tensor product of $X_1,\ldots,X_n$
will be denoted by  $\tensor{X_1}{X_n}$ { and  the space  of multilinear bounded operators from $\cartesian{X}{1}{n}$ to $Y$
     by $ \mathcal{L}(\spaces{X}{1}{n}; Y)$. When $Y$ is the scalar field, we will use the notation $\mathcal{L}(\spaces{X}{1}{n}).$}
The general theory of Banach spaces that we will  use  can be found in
\cite{DJT} and \cite{Ryan-libro}.

We  briefly recall now
   how $\Sigma$-operators are used  to study bounded multilinear mappings. Details can be found in \cite{MFU}.  Each  multilinear bounded operator $T\in \mathcal{L}(\spaces{X}{1}{n};Y)$ admits a  factorization $T=\hat{T}\circ \otimes$, where  $\hat{T} \in  \mathcal{L}(\Tensor{X}{1}{n}{};  Y)$ is   such     that  for every $x_i\in X_i \; i=1,\ldots, n $, \;
    $T(\spaces{x}{1}{n})=\hat{T}(\Tensor{x}{1}{n}{}).$
  The mapping $f_T:=\hat{T}_{|_{\Sigma}}$, where ${\Sigma}_{X_1,\ldots,X_n}:= \{\Tensor{x}{1}{n}{}\in  \Tensor{X}{1}{n}{}; \; x_i\in X_i\}$ is the   Segre cone of the spaces $X_i$, is called a {\sl $\Sigma$-operator}.  The space of continuous $\Sigma$-operators endowed with the Lipschitz norm will be denoted $\mathcal{L}\left(\Sigma_{\spaces{X}{1}{n}}; Y\right)$ and $\mathcal{L}\left(\Sigma_{\spaces{X}{1}{n}}\right)$ in the case where $Y$ is the scalar field.
The following mappings are  isometric isomorphisms (see \cite[Theorem 3.2]{MFU}):
 \[
 \begin{array}{ccccc}
  \mathcal{L}\left({\spaces{X}{1}{n}}; Y \right) & \stackrel{\Psi}{\longrightarrow} & \mathcal{L}\left(\Sigma_{\spaces{X}{1}{n}}; Y\right) &\stackrel{\Phi}{\longrightarrow} & \mathcal{L}\left({\CTensor{X}{1}{n}{\pi}}; Y \right) \\
   T & \mapsto &  f_T & \mapsto   & \hat{T}
 \end{array}
 \]

The general procedure to  go  from a given  theory on linear operators  to the broader context of  $\Sigma$-operators  is done in two steps: the   first step is to interpret a specific type of   boundedness condition   on linear operators $\{S : X\rightarrow  Y \}$ as a  continuous condition, thus, as a Lipschitz condition. The second step is to  formulate  such Lipschitz   condition for  $\Sigma$-operators $\{f : \Sigma_{\spaces{X}{1}{n}} \rightarrow  Y \}$.
In the case of the $p$-summability property, this procedure gives rise to the following definition:

$f\in\mathcal{L}(\Sigma_{X_1,\ldots,X_m}; Y)$  is  an  {\sl absolutely $p$-summing} $\Sigma$-operator    if there is a
constant $c\geq 0$ such that for every   $i=1,\ldots,k$
and  every  $u_i, v_i \in \Sigma_{\spaces{X}{1}{n}}$, the following inequality holds:
\begin{multline}\label{eq: def sigma p sumante}
   {\left(\sum_{i=1}^k \left\|{f}(u_i)-{f}(v_i)\right\|^p \right)^{1/p}  }
 \leq  \\  c\cdot \sup \left\lbrace \left (\sum_{i=1}^k \left|\varphi(u_i)- \varphi(v_i)\right|^p\right)^{1/p}; \; \varphi\in B_{\mathcal{L}\left(\Sigma_{\spaces{X}{1}{n}}\right)}\right\rbrace
\end{multline}

 When $n=1$,  we have that $\Sigma_X=X$.  In this case  absolutely $p$-summing $\Sigma$-operators and   absolutely $p$-summing  linear operators coincide, as well as
  $\Sigma$-operators and  bounded linear mappings coincide.

 The space of absolutely $p$-summing $\Sigma$-operators, denoted by   $\Pi_{p}(\Sigma_{\spaces{X}{1}{n}}; Y)$,  is a Banach space when defining  the norm $\pi_{p}(f)$ as the smallest constant $c \geq 0$ satisfying  inequality (\ref{eq: def sigma p sumante}).

The following result will  be used in  Theorem \ref{thm: equiv p summ multi op}. It can be directly proved.
\begin{lemma}\label{lem: Sigma-inclusion} Let $\delta_{\Tensor{x}{1}{n}{}}$  be defined as $\delta_{\Tensor{x}{1}{n}{}}(\varphi):=\varphi(\spaces{x}{1}{n})$ for every $\varphi\in \mathcal{L}({X_1,\ldots, X_n})$. Then,  the following map is an isometric inclusion.
\[
 \begin{array}{ccccc}
\Sigma_{X_1,\ldots, X_n} &\stackrel{i}{\hookrightarrow} &  C(B_{\mathcal{L}(X_1, \ldots, X_n)}),w^*)
   \\
  {\Tensor{x}{1}{n}{}}& \mapsto  & \delta_{\Tensor{x}{1}{n}{}}.
 \end{array}
 \]
\end{lemma}

\section{Lipschitz $p$-summing multilinear operators.}\label{sec: Lips multi}

We will denote   $j_p: C(B_{\mathcal{L}\left(X_1,\ldots,X_n\right)},w^*)\rightarrow   L_p\left(\mu\right)$ and  $i_p: L_{\infty}(\mu)\rightarrow   L_p\left(\mu\right)$ the natural inclusion mappings and
$i_Y: Y\rightarrow \ell_{\infty}^{B_Y^*}$ the natural isometric
inclusion $i_Y(z)=(x^*(z))_{x^*\in B_Y^*}$.

Now we apply  the  geometrical approach explained in Section \ref{sec: sigma-op} to obtain the following fundamental equivalences:

\begin{theorem}\label{thm: equiv p summ multi op} Let  $1 \leq p < \infty$ and  $T:X_1\times\cdots\times X_n\rightarrow Y$ be an $n$-linear operator between Banach spaces.       The  following conditions for $T$ are equivalent:
\begin{enumerate}
\item[i)]  There exists $c>0$ such that for $k\in\mathbb{N}$, $i=1,\ldots,k$ and $u_i, v_i \in \cartesian{X}{1}{n}$,
  \begin{multline*} \label{multi defsumantes.1}
  {\left(\sum_{i=1}^k \left\|{T}\left(u_i\right)-T \left(v_i\right)\right\|^p \right)^{1/p}  }
  \leq  \\  c\cdot \sup \left\lbrace \left(\sum_{i=1}^k \left|{{\varphi}}\left(u_i\right)-{\varphi}\left(v_i\right)\right|^p\right)^{1/p}; \; \varphi \in B_{\mathcal{L}\left(\spaces{X}{1}{n}\right)}\right\rbrace.
 \end{multline*}

\item[ii)] There is a constant $c>0$ and a regular probability measure $\mu$ on the space $\left(B_{\mathcal{L}\left(X_1,\ldots,X_n\right)},w^*\right)$ such that for each $u=(u_1,\ldots,u_n), v=(v_1,\dots, v_n) \in \cartesian{X}{1}{n} $ we have that
\begin{equation*} \label{eqn: domination}
{\left\|T\left(u\right)-T\left(v\right)\right\| }
 \leq  c\cdot\left(\int_{B_{\mathcal{L}\left(X_1,\ldots,X_n\right)}} \hspace{-1.7cm}\left|\varphi\left(u\right)-\varphi\left(v\right)\right|^p\,d\mu(\varphi)\right)_{\ .}^{1/p}
\end{equation*}

 \item[iii)]  There exist a regular Borel probability measure $\mu$ on  $(B_{\mathcal{L}\left(X_1,\ldots,X_n\right)},w^*)$,
 a  subset $\Sigma_p:=\left(j_p\circ \otimes\right)\left({X_1,\ldots, X_n}\right)$ of $L_p\left(\mu\right)$  and a Lipschitz function ${{h}_T}: {\Sigma_p}\rightarrow {Y}$ such that
  $T=h_{T} \circ j_p \circ \otimes$, that is, in such a way that the following diagram commutes:
   \[
\begin{tikzcd}\label{p-summ diag}
        {\cartesian{X}{1}{n}} \arrow{r}{T}\arrow{d}{\otimes} &  Y\\
        \Sigma_{X_1,\ldots, X_n}  \arrow{r}{{j_p}_{|_{\Sigma}}}\arrow[phantom,"\cap"]{d} &  \Sigma_p \arrow[phantom,"\cap"]{d}\arrow{u}[swap]{h_T}\\
        C(B_{\mathcal{L}\left(X_1,\ldots,X_n\right)},w^*)\arrow{r}{{j_p}}   &    L_p\left(\mu\right).\\
 \end{tikzcd}
\]

\item[iv)]
There  exist  a  probability space $(\Omega,\Sigma,\mu)$, a  multilinear
operator $\nu: X_1\times\ldots \times X_n\rightarrow L_{\infty}(\mu)$, $\|\nu\|=1$,  and   a Lipschitz
function $\funcion{\tilde{h_T}}{L_p\left(\mu\right)}{\ell_{\infty}^{B_{Y^*}}}$
such that the following diagram commutes:

     \[
\begin{tikzcd}\label{p-summ multi  linfty}
        X_1\times\ldots \times X_n \arrow{r}{T}\arrow{dd}[swap]{\nu} &  Y \arrow{rd}{i_Y} & \\
        & & \ell_{\infty}^{B_{Y^*}} \\
        L_{\infty}(\mu)\arrow{r}{{i_p}}   &    L_p\left(\mu\right)\arrow{ru}{\tilde{h_T}}. & \\
 \end{tikzcd}
\]


\end{enumerate}

 If $\pi_p^{Lip}(T):=\inf\{c;\; i) \; \mbox{holds}\}$,
   then  $\pi_p^{Lip}(T)  =\inf\{c;\; ii) \; \mbox{holds}\}$ and   $\normasub{\tilde{h_T}}{Lip}= \pi_p^{Lip}(T)$ when the spaces are  real  and  $\pi_p^{Lip}({T}) \leq \normasub{\tilde{h_T}}{Lip}\leq \sqrt{2}\pi_p^{Lip}(T)$ when the spaces are complex.

\end{theorem}

\begin{pf}
The  proof is  like the proof of the original (linear) Pietsch factorization Theorem.
 {
Indeed,  for proving  the equivalence between $(i)$ and $(ii)$, it is possible to apply \cite[Theorem 2.2]{BotPelRue2}. There the authors extend the original statement concerning bounded linear operators (see \cite[Theorem 2.12]{DJT}) to a broader class of mappings.
 }

To  see
  $(ii) \Rightarrow (iii)$, we will use Lemma \ref{lem: Sigma-inclusion} without mentioning it.  Consider a   regular Borel probability measure  $\mu$  as in $(ii)$.  Let
 $\Sigma_p:=\left(j_p\circ \otimes\right)\left(\cartesian{X}{1}{n}\right) \subset L_p\left(\mu\right)$ and  define    $h_T : \Sigma_p \longrightarrow Y $
as $h_T((j_p\circ \otimes)(\spaces{u}{1}{n})):=T(\spaces{u}{1}{n})$ for every $(\spaces{u}{1}{n})\in {X_1\times\ldots\times X_n}$.    $h_T$ is well defined, since whenever $(j_p\circ \otimes)(\spaces{u}{1}{n})=(j_p\circ \otimes)(\spaces{v}{1}{n})$, $(ii)$ guarantees that
\begin{multline*}\|T(\spaces{u}{1}{n})-T(\spaces{v}{1}{n})\|\leq \\ \|(j_p\circ \otimes)(\spaces{u}{1}{n})-(j_p\circ \otimes)(\spaces{v}{1}{n})\|_{L_p(\mu)}=0. \end{multline*}
Finally,  $\|{h_T}\|_{{Lip}} = \pi_p^{Lip}(T)$ holds because for every $z,w \in \Sigma_p$,
$${\left\|h_T\left(z\right)-h_T\left(w\right)\right\| }
 \leq  \|h_T\|_{Lip}\cdot\left(\int_{B_{\mathcal{L}}} \left|\varphi\left(z\right)-\varphi\left(w\right)\right|^p\,d\mu(\varphi)\right)_{\ .}^{1/p}$$
 and  $\pi_p^{Lip}(T)$ is the infimum among the constants satisfying $(i)$.

$(iii) \Rightarrow (iv)$. Assume that ${{h}_T}: {\Sigma_p}\rightarrow {Y}$ is as in $(iii)$.  Then,  ${i_Y \circ h_{T}}$ is a Lipschitz function with $\normasub{i_Y \circ h_{T}}{Lip}=\normasub{h_{T}}{Lip}=\pi_p^{Lip}(T)$. { Using the Lipschitz 1-injectivity property of $\ell_{\infty}^\Gamma$ (see  \cite[Lemma~1.1]{BenyLind})},  we  find a Lipschitz extension $\tilde{h_T}$ of $i_Y \circ {h_T}$ defined on  $L_p(\mu)$,  such that
$
\normasub{\tilde{h_T}}{Lip}=\normasub{i_Y \circ h_{T}}{Lip}=\pi_p^{Lip}(T)
$
when the spaces were  real spaces,  and
$
\pi_p^{Lip}(T)\leq \normasub{\tilde{h_T}}{Lip}\leq
\sqrt{2}\cdot\pi_p^{Lip}(T)
$
when the spaces were complex.

$(iv) \Rightarrow (ii)$. For  $(\spaces{u}{1}{n}),(\spaces{v}{1}{n})\in \cartesian{X}{1}{n}$,  inequality  ${ (ii)}$, and consequently  $(i)$, follows from:
\begin{multline*}
\|T (\spaces{u}{1}{n}) -T(\spaces{v}{1}{n})\|=\|(i_Y\circ T) (\spaces{u}{1}{n}) -(i_Y \circ T)(\spaces{v}{1}{n})\| \\
=\|\tilde{h_T}\circ i_p \circ \nu(\spaces{u}{1}{n}) -\tilde{h_T}\circ i_p \circ  \nu(\spaces{v}{1}{n})\|
\\\leq \|\tilde{h_T}\|_{Lip}\|i_p \circ  \nu(\spaces{u}{1}{n}) - i_p \circ  \nu(\spaces{v}{1}{n})\|_{L_p(\mu)}.
\end{multline*}
 \end{pf}

  \begin{definition}\label{def p summing multi} A multilinear operator $T\in\mathcal{L}(X_1,\ldots,X_n; Y)$ is  {\bf Lipschitz $p$-summing}  if it satisfies any of the equivalent conditions in Theorem \ref{thm: equiv p summ multi op}. The Lipschitz $p$-summing norm of $T$ is defined as $\pi_p^{Lip}(T)$.
\end{definition}

$\Pi_{p}^{Lip}(\spaces{X}{1}{n};Y)$ will denote       the Banach space  of  Lipschitz $p$-summing multilinear operators with  the norm $\pi_p^{Lip}$.   Clearly,  $\Pi_{p}^{Lip}(\spaces{X}{1}{n};Y)$ and $\Pi_p(\Sigma_{\spaces{X}{1}{n}};Y)$ are isometrically isomorphic  Banach spaces, via the mapping $T\mapsto f_T$.

  \begin{remark}  We use the name  {\sl Lipschitz} in Definition \ref{def p summing multi} to   highlight the role that the differences play in such notion. We also try   to avoid  confusion with other notions   of $p$-summability of multilinear operators that appear in the literature. The main reason is, however, that $\Sigma$-operators are Lipschitz functions, although non-zero multilinear mappings are not.  In Subsection \ref{subs: Lips} we discuss this fact.
\end{remark}

\begin{remark}\label{rmk: factoriz}
Other equivalent characterizations of Lipschitz $p$-summing operators can be proved.
We mention one that will be used later.
Its  proof relies on $(iii)$ and  {  the  Lipschitz $1$-injectivity property of  $\ell_{\infty}^{\Gamma}$: }

There  exist  a  probability space $(\Omega,\Sigma,\mu)$ and   a Lipschitz
function ${\tilde{h_T}}$
such that the following diagram commutes:

     \[
\begin{tikzcd}\label{p-summ multi  linfty2}
        X_1\times\ldots \times X_n \arrow{r}{T}\arrow{dd}[swap]{i\circ \otimes} &  Y \arrow{rd}{i_Y} & \\
        & & \ell_{\infty}^{B_{Y^*}} \\
        C(B_{\mathcal{L}\left(X_1,\ldots,X_n\right)},w^*)\arrow{r}{{j_p}}   &    L_p\left(\mu\right)\arrow{ru}{\tilde{h_T}}. & \\
 \end{tikzcd}
\]
\end{remark}

\subsection*{Examples of  Lipschitz $p$-summing operators.}
  Clearly, any bounded multilinear form is Lipschitz $p$-summing.
  This implies that any multilinear operator  whose image belongs to a
   finite dimensional space is Lipschitz $p$-summing.
   Other natural examples   are the following:

 Let $K$, $K_1,\ldots, K_n$  be  compact Hausdorff spaces and let   $\mu$, $\nu$  be  positive regular Borel measures      on $K$ and  $K_1\times\cdots\times K_n$  respectively. If   $h\in L_p(\mu)$ and  $f\in L_p(\nu)$, then the  multilinear operators  $T_h\in\mathcal{L}({C(K),\ldots, C(K)}; L_p(\mu))$, and   $ S_f\in\mathcal{L}({C(K_1), \ldots,  C(K_n)}; L_p\left(\nu\right))$  defined as $T_h(g_1, \ldots,  g_n):=h(w)\cdot g_1(w)\cdots g_n(w)$   and
 $S_f(g_1, \ldots,  g_n):=f(w_1,\ldots,w_n)\cdot g_1(w_1)\cdots g_n(w_n)$, are  Lipschitz $p$-summing.

Given  $\lambda=(\lambda_k)_{k}\in \ell_p$, the
 diagonal operator
  $T_{\lambda}\in\mathcal{L}({\ell_{\infty}, \ldots, \ell_{\infty}};  \ell_p)$, defined as
$   T_{\lambda}\left((a_k^1)_{k},\ldots,({a_k^n})_{k}\right):=({\lambda_k\cdot a_k^1\cdots a_k^n})_{k}$, is Lipschitz $p$-summing.

\subsection*{Main properties.} To  establish  the generalizations of some of the main results in the linear case,  first we see that
  Lipschitz $p$-summing operators are located between the following    two  multilinear  generalizations  of absolutely $p$-summing linear operators. $\mathcal{L}^p_{ss}$ denotes the set of strongly $p$-summing multilinear operators (see its definition in  \cite{Dimant03}):
{ \begin{proposition}\label{prop: lineal psumm es Lipsch psumm}  For $X_1,\ldots,X_n, Y $ Banach spaces,
 \begin{equation*}\label{eq: strongly}
 \begin{split}
  \{T\in  \mathcal{L}(\spaces{X}{1}{n}&; Y); \; \hat{T}\in \Pi_p(\CTensor{X}{1}{n}{\pi} ; Y)\}  \subset \\ & \Pi_{p}^{Lip}(\spaces{X}{1}{n}; Y)\;\subset\;\mathcal{L}^p_{ss}(\spaces{X}{1}{n};Y).
     \end{split}
      \end{equation*}
\end{proposition}}
\begin{pf}
{ To prove the first inclusion, let us consider a bounded multilinear mapping $T$ such that its associated linear operator $\hat{T}$ defined on $\CTensor{X}{1}{n}{\pi}$ is an absolutely $p$-summing operator. Then, the for any finite set $z_1,\ldots, z_k\in \CTensor{X}{1}{n}{\pi}$ it holds that
 \begin{equation*} \label{multi defsumantes.1}
  {\left(\sum_{i=1}^k \left\|\hat{T}\left(z_i\right)\right\|^p \right)^{1/p}  }
  \leq    \pi_p(\hat{T})  \sup \left\lbrace \left(\sum_{i=1}^k \left|{{\varphi}}\left(z_i\right)\right|^p\right)^{1/p}; \; \varphi \in B_{\mathcal{L}\left(\CTensor{X}{1}{n}{\pi}\right)}\right\rbrace.
 \end{equation*}
  So, in particular, this estimate holds for any finite set of vectors of the form
   $z_i= u^i_1\otimes \ldots \otimes u^i_n-v^i_1\otimes \ldots \otimes v^i_n\in \Sigma_{\spaces{X}{1}{n}}-\Sigma_{\spaces{X}{1}{n}}$.}
 This  gives rise to $(i)$ in Theorem \ref{thm: equiv p summ multi op}. To see the other relation,  just observe  that  $T\in \mathcal{L}(\spaces{X}{1}{n}; Y)$ is  {\sl strongly $p$-summing} exactly when $T$ satisfies    condition $(i)$ in Theorem \ref{thm: equiv p summ multi op}  for every $k\in\mathbb{N}$, $i=1,\ldots,k$,  $u_i \in \cartesian{X}{1}{n}$ and   $v_i=0$.
\end{pf}

Example 3.3 in  \cite{CarDimant03}  shows that $T\in\mathcal{L}(\ell_1,\ell_1; \ell_1)$,  $T((x),(y))=(x_jy_j)_j$ is   a strongly $p$-summing multilinear operator whose associated linear mapping is not absolutely $p$-summing.  We do not know if $T$ is Lipschitz $p$-summing.

\begin{corollary}\label{coro: Aron-Berner} The  Aron-Berner extension of a Lipschitz $p$-summing multilinear operator $T\in\Pi_{p}^{Lip}(\spaces{X}{1}{n}; Y)$  takes its values in $Y$.
\end{corollary}
\begin{pf}
 By Proposition \ref{prop: lineal psumm es Lipsch psumm}, $T$ is strongly $p$-summing. Thus, by  Theorem 2.2 in \cite{Dimant03}, its Aron-Berner extension is in $\mathcal{L}(\spaces{X^{**}}{1}{n}; Y)$.
\end{pf}
The  statements  stated below as propositions  can be proved using  standard arguments.
We will only  prove their corollaries. This type of results are usually  known
 as {\sl coincidence theorems}.

\begin{corollary}[A Grothendieck Theorem for  multilinear operators]\label{coro: Groth Thm} If $\spaces{X}{1}{n}$, $n\in\N$,  are $\mathcal{L}_1$ spaces and $Y$ is an $\mathcal{L}_2$ space, then
   $$\Pi_{1}^{Lip}(\spaces{X}{1}{n};Y)=\mathcal{L}(X_1,\ldots,X_n; Y).$$
  \end{corollary}
\begin{pf}

 Under our assumptions, the space $\CTensor{X}{1}{n}{\pi}$ is also an $\mathcal{L}_1$ space (see \cite[p.314]{DeFlo}). By  Grothendieck's Theorem  \cite[Theorem 3.1]{DJT} and Proposition \ref{prop: lineal psumm es Lipsch psumm},   this implies   that $$\{T; \hat{T}\in\Pi_{1}(\CTensor{X}{1}{n}{\pi};Y)\}=\Pi_{1}^{Lip}(\spaces{X}{1}{n};Y)=\mathcal{L}(X_1,\ldots,X_n; Y).$$
\end{pf}

\begin{proposition}[Inclusion Theorem]\label{prop: inclusion thm}
 \label{prop: multi usual p,q contentions}
If $1 \leq p \leq q < \infty$, then $$\Pi_{p}^{Lip}(\spaces{X}{1}{n};Y) \subset
\Pi_{q}^{Lip}(\spaces{X}{1}{n};Y) $$  and $\pi_{q}^{Lip}(T)\leq\pi_{p}^{Lip}(T)$\ for every $T \in
\Pi_{p}^{Lip}(\spaces{X}{1}{n};Y)$.
\end{proposition}

The ideal properties of absolutely $p$-summing linear operators extend to Lipschitz $p$-summing operators as follows:

\begin{proposition}\label{prop: multi-ideal}  Let $1\leq p <\infty$. For every $n\in\N$ and
every $X_1,\ldots,X_n, Y,Z$ Banach spaces, it holds that
\begin{enumerate}
\item[i)] $\Pi_{p}^{Lip}(\spaces{X}{1}{n};Y)$ is a linear subspace of  $\mathcal{L}(\spaces{X}{1}{n};Y)$ and $\|T\| \leq \pi_{p}^{Lip}(T)$.
\item[ii)] The space
 $(\Pi_{p}^{Lip}(\spaces{X}{1}{n};Y), \pi^{Lip}_{p} )$ is a Banach space.
 It contains the multilinear operators  whose range lies  in a finite dimensional subspace of $Y$.

\item[iii)]  Let $T \in \Pi^{Lip}_{p}({X_1,\ldots,X_n};Y)$, $R \in \mathcal{L}\left(Y,Z\right)$ and $S_i \in \mathcal{L}\left(W_i,X_i\right),$
 $i=1\ldots,n$. Then $R \circ T \circ(\cartesian{S}{1}{n})\in \Pi^{Lip}_{p}(\spaces{W}{1}{n};Z)$ and
 $\pi^{Lip}_{p}(R \circ T \circ(\cartesian{S}{1}{n}) )\leq \norma{R}\cdot\pi^{Lip}_{p}(T)\cdot\norma{S_1}\cdots\norma{S_n}$.

\item[iv)] $\pi^{Lip}_{p}(\Lambda_n)=1$,  where ${\Lambda_n: \cartesian{\mathbb{K}}{}{}} \rightarrow \mathbb{K}$  is  $\Lambda_n((z_1,\cdots ,z_n))=z_1\cdots z_n$.

\item[v)] Let $T \in \mathcal{L}({X_1,\ldots,X_n};Y)$ and $R\in\Pi_{p} \left(Y,Z\right)$.
 Then it holds that  $R\circ T \in \Pi^{Lip}_{p}({X_1,\ldots,X_n};Z)$ and $\pi^{Lip}_{p}(R \circ T)\leq
 \pi_{p}(R)\cdot \|T\|_{Lip}$.

 \item[vi)] If $\funcion{i}{Y}{Z}$ is a linear isometry, then   $T\in \mathcal{L}({X_1,\ldots,X_n};Y)$
 is Lipschitz $p$-summing if and only if $i \circ T \in \mathcal{L}({X_1,\ldots,X_n}; Z)$ is Lipschitz $p$-summing
  and $\pi^{Lip}_{p}(T)=\pi^{Lip}_{p}(i \circ T)$.
\end{enumerate}
\end{proposition}

\begin{corollary}\label{coro pi1 = pi2} If $\spaces{X}{1}{n}$, $n\in\N$,  are $\mathcal{L}_1$ spaces, then
   $$\Pi_{1}^{Lip}(\spaces{X}{1}{n};Y)=\Pi_{2}^{Lip}(\spaces{X}{1}{n};Y).$$
  \end{corollary}
\begin{pf} By Proposition \ref{prop: inclusion thm}, we only need to check that $\Pi_{2}^{Lip}(\spaces{X}{1}{n};Y)\subset\Pi_{1}^{Lip}(\spaces{X}{1}{n};Y).$
Let $T\in \Pi_{2}^{Lip}(\spaces{X}{1}{n};Y)$ and let $i_Y\circ T= \tilde{h_T}\circ j_2 \circ (i\circ \otimes)$ be the factorization given in Remark \ref{rmk: factoriz}. By Corollary \ref{coro: Groth Thm},
 $j_2 \circ (i\circ \otimes)$  is Lipschitz $1$-summing.  This implies that $\tilde{h}_T\circ j_2 \circ (i\circ \otimes)$, which  equals $i\circ T$,  is Lipschitz $1$-summing. By $(vi)$ Proposition \ref{prop: multi-ideal}, $T$ is Lipschitz $1$-summing, too.
 \end{pf}

 We  will use the notation $x_1\otimes \stackrel{{j}}{\stackrel{\vee}\cdots}\otimes x_n$  meaning $x_1\otimes\cdots \otimes x_{j-1}\otimes  x_{j+1}\otimes\cdots\otimes x_n${, $X_1,\stackrel{{i_1,\ldots, i_l}}{\stackrel{\vee}\ldots},  X_n$ meaning the product of the spaces $X_i$, $i\notin \{i_1,\ldots, i_l\}$   and   $T_{\spaces{x^0}{i_1}{i_l}}\in \mathcal{L}(X_1,\stackrel{{i_1,\ldots, i_l}}{\stackrel{\vee}\ldots},  X_n)$  to denote the mapping such that if  each $x^0_{i_j}$ is placed in the $i_j$-th position, it coincides whith  $T$.}

\begin{proposition}\label{prop: less degree}
Let  $1 \leq p< \infty$ and $x^0_{i_j}\in X_{i_j}, j=1,\dots,l$, $l<n$  be fixed. If $T\in \Pi_{p}^{Lip}(\spaces{X}{1}{n};Y)$  then, $T_{\spaces{x^0}{i_1}{i_l}}\in  \Pi_{p}^{Lip}(X_1,\stackrel{{i_1,\ldots, i_l}}{\stackrel{\vee}\ldots}, X_n;Y)$, $\pi_p^{Lip}(T_{\spaces{x^0}{i_1}{i_l}})\leq\|x^0_{i_1}\|\cdots \|x^0_{i_l}\|\pi_p^{Lip}(T)$.
\end{proposition}

\begin{corollary}[A weak Dvoretzky-Rogers Theorem ]\label{coro: teorema debil de Dvoretzky-Rogers} Let  $1 \leq p< \infty$ and let $X$ be a Banach space. Then
$\Pi_{p}^{Lip}(X,\stackrel{n}{\ldots}, X; X) =\mathcal{L}\left(^n X;X\right)$ if and only if
 $X$ is a finite dimensional space.
\end{corollary}
\begin{pf} Let us assume that every multilinear operator is Lipschitz $p$-summing. Whenever
 $l=n-1$,   Proposition \ref{prop: less degree}  guarantees  that any  linear mapping of the form $T_{\spaces{x^0}{i_1}{i_l}}$ is absolutely $p$-summing. By the weak Dvoretzky-Rogers  theorem, this already implies that $X$ is finite dimensional.
 The reverse implication follows from  $ii)$ in Proposition \ref{prop: multi-ideal}.
\end{pf}

\begin{corollary}[A Lindenstrauss-Pe\l czy\'nski Theorem]\label{coro l1 Hilbert} Let  $\spaces{X}{1}{n}$, $Y$ be  Banach spaces, each  $\spaces{X}{1}{n}$  with an unconditional basis. Each   $X_i$ is isomorphic to  $\ell_1^{{\Gamma}_i}$ for some  $\Gamma_i$ and  $Y$ is isomorphic to a Hilbert space   if and only if
 $$\Pi_{1}^{Lip}(\spaces{X}{1}{n}; Y) = \mathcal{L}\left(\spaces{X}{1}{n};Y\right).$$
\end{corollary}
 \begin{pf}
  If    $X_i \simeq \ell_1^{\Gamma_i}$ and  $Y\simeq H$,  the result is given by   Corollary \ref{coro: Groth Thm}.
Assume now that each  $\spaces{X}{1}{n}$ has an   unconditional basis and, together  with $Y$  are such that every $n$-linear operator is Lipschitz $1$-summing. Let us fix $n$ norm one vectors  $x_i^0\in X_i$  and $n$ norm one linear projections $P_i\in\mathcal{L}(X_i,X_i)$ onto $\langle {x_i^0}\rangle, $ respectively.   Given    an arbitrary  $S\in \mathcal{L}(X_1; Y)$,  let us define ${S^1}  \in \mathcal{L}(\spaces{X}{1}{n}; Y)$ as
$S^1 (\spaces{x}{1}{n})=S(x_1)P_2(x_2)\cdots P_n(x_n)$. Under our assumption, $S^1 $ is  Lipschitz $1$-summing. Then, by Proposition \ref{prop: less degree}, we have that $S^1_{x_2^0,\dots, x_n^0}$ is an absolutely $1$-summing linear operator. Finally, observe that  $S^1_{{x_2^0,\dots, x_n^0}}=S$. We have proved that   $\mathcal{L}(X_1; Y)=\Pi_1(X_1,Y).$  By  \cite[Theorem 4.2]{LP}, we know that  $X_1$ is isomorphic to  an $\ell_1^{\Gamma_1}$ and $Y$ is isomorphic to a Hilbert space. An analogous argument can be done with the other indexes $i=2,\ldots,  n$.
\end{pf}

\section[Dunford-Pettis and Lipschitz $p$-summing]{Dunford-Pettis and Lipschitz $p$-summing multilinear operators }\label{sec: DP}

A multilinear mapping is said to be  {\sl compact} if it maps bounded sets into relatively compact sets.
The  space of compact multilinear  operators  will be denoted $ \mathcal{K}(X_1,\ldots,X_m; Y).$  It holds that
 ${T}\in \mathcal{K}(X_1,\ldots,X_n; Y) $ if and only if $\hat{T}\in\mathcal{K}(\CTensor{X}{1}{n}{\pi}; Y). $ Let us see that  in general,   there is not a  containment relationship  between compact and Lipschitz $p$-summing multilinear operators:

\begin{example} Consider  $S\in  {\Pi_p}(\ell_1; \ell_2)\setminus\mathcal{K}(\ell_1; \ell_2)$;  for instance, the formal inclusion mapping.  Let ${S^1}  \in \mathcal{L}(\ell_1,X_2,\ldots,X_n; \ell_2)$ be the  multilinear mapping  $S^1 (\spaces{x}{1}{n})=S(x_1)P_2(x_2)\cdots P_n(x_n)$,  constructed from $S$ as in the proof of Corollary \ref{coro l1 Hilbert}.
 $S^1$ is  a non compact Lipschitz $p$-summing multilinear mapping.
\end{example}

\begin{example} Let  $S\in  \mathcal{K}(\ell_2; \ell_2)\setminus{\Pi_p}(\ell_2; \ell_2)$ (see an example in \cite[p.38]{DJT}). Constructing $S^1$ as before, $S^1 \in \mathcal{K}(\ell_2,X_2,\ldots,X_n; \ell_2)\setminus \Pi_p^{Lip}(\ell_2,X_2,\ldots,X_n; \ell_2).$
\end{example}

  Proposition \ref{prop: psum DP} states  that, however,  Lipschitz $p$-summing operators are of the following type:

\begin{definition}\label{def: DP} We say that a multilinear operator  $T\in\mathcal{L}(X_1,\ldots,X_m; Y)$ is { Dunford-Pettis} if  its associated $\Sigma$-operator $f_T$ transforms the weakly Cauchy sequences in $\Sigma_{X_1,\ldots,X_m}$ into  norm convergent sequences in $Y$. In this case, we say that the $\Sigma$-operator $f_T$ is Dunford-Pettis.   The space of Dunford-Pettis multilinear operators will be denoted  $\mathcal{DP}(X_1,\ldots,X_m; Y)$.
 \end{definition}

 In the case of linear mappings, $n=1$, Definition \ref{def: DP} coincides with the usual definition  of a Dunford-Pettis (or completely continuous) linear operator (see \cite{Bombal90} and the references therein). A  notion of Dunford-Pettis for  polynomials was   studied in \cite{MFU2001}.

\begin{proposition}\label{prop: linear-multi DunfPet}  Let $\spaces{X}{1}{n}$ be Banach spaces and $1\leq p <\infty$. Then,
\begin{enumerate}
  \item $ \mathcal{K}(X_1,\ldots,X_n; Y)\subset \mathcal{DP}(X_1,\ldots,X_n; Y)$.

  \item If $T\in\mathcal{L}(X_1,\ldots,X_n; Y)$  is such that $\hat{T}\in \mathcal{DP}(\CTensor{X}{1}{n}{\pi}; Y)$, then $T \in \mathcal{DP}(X_1,\ldots,X_n; Y)$.
\end{enumerate}
\end{proposition}
\begin{pf} $\mathit{(1)}$ follows from a compactness argument. To prove $\mathit{(2)}$, let
$\hat{T}$  be such that it transforms the  weakly Cauchy sequences of  $\CTensor{X}{1}{n}{\pi}$ into norm convergent sequences in $Y$.   In particular, this holds for  the weakly Cauchy sequences in $\Sigma_{\spaces{X}{1}{n}}$. Since  $f_T=\hat{T}_{|_{\Sigma_{\spaces{X}{1}{n}}}}$, we have the desired condition on  $f_T$.
\end{pf}

The following theorem generalizes an important and non-trivial result, which corresponds to the case $n=1$ of our result. We denote $c_0$ the space of sequences converging to zero with the norm $\|\cdot\|_{\infty}$.

 \begin{theorem}\label{thm: DP}  Let $\spaces{X}{1}{n}$ be Banach spaces. The following are equivalent:
 \begin{enumerate}
   \item  The space $\CTensor{X}{1}{n}{\pi}$ contains no isomorphic copies of $\ell_1$.
   \item Every bounded sequence in $\Sigma_{\spaces{X}{1}{n}}$ contains a weakly Cauchy subsequence.
   \item  For every Banach space $Y$,   $\mathcal{K}(X_1,\ldots,X_n; Y)=\mathcal{DP}(X_1,\ldots,X_n; Y).$
    \item $ \mathcal{K}(X_1,\ldots,X_n; c_0)=\mathcal{DP}(X_1,\ldots,X_n; c_0).$
 \end{enumerate}
 \end{theorem}
 \begin{pf}
By  Rosenthal's dichotomy Theorem,   $\mathit{(1)}$ is equivalent to the fact that every bounded sequence in $\CTensor{X}{1}{n}{\pi}$ contains a weakly Cauchy subsequence. In particular, this holds for the subset $\Sigma_{\spaces{X}{1}{n}}$, which gives $\mathit{(2)}$. Assume $\mathit{(2)}$.  By Proposition \ref{prop: linear-multi DunfPet}, to prove $\mathit{(3)}$ it is enough to prove $\supset$. Let $T$ be a Dunford-Pettis multilinear operator and
{let $A_i\in X_i, \; i=1,\ldots, n$ be arbitrary bounded subsets. Consider the   subsets $A:=\cartesian{A}{1}{n}\subset \cartesian{X}{1}{n}$ and  $A^{\otimes}:=\{\Tensor{x}{1}{n}{}; \; x_i\in A_i\}\subset \Sigma_{\spaces{X}{1}{n}}$.}
 By $\mathit{(2)}$, any sequence in $A^{\otimes}$ has  a weakly Cauchy subsequence $\{\Tensor{a^j}{1}{n}{}\}_j\subset A^{\otimes}$. Since  $T$ is  Dunford-Pettis,  $\{f_T(\Tensor{a^j}{1}{n}{})\}_j$ is  a norm convergent sequence. But $\{f_T(\Tensor{a^j}{1}{n}{})\}_j=\{T(\spaces{a^j}{1}{n})\}_j$. Consequently, $T$ is compact. $\mathit{(4)}$ is a particular case of $\mathit{(3)}$. Finally, let us  see that  $\mathit{(4)}$ implies $\mathit{(1)}$.    If   $\ell_1\subset   \CTensor{X}{1}{n}{\pi}$,  there exist  a  linear operator $\hat{S} \in\mathcal{DP}(\CTensor{X}{1}{n}{\pi}; c_0)\setminus  \mathcal{K}(\CTensor{X}{1}{n}{\pi}; c_0)$ (for a proof, see Proposition \cite[II.6]{Bombal90}). Proposition \ref{prop: linear-multi DunfPet} implies that its associated multilinear  operator   $S$ is Dunford-Pettis. Since  $S$ is non compact, $\mathit{(4)}$ does not hold.
 \end{pf}

 The linear case $n=1$ of the following results can be found in \cite[Corollary 1.7]{Pisier}.
\begin{proposition}\label{prop: psum DP} Let $1\leq p < \infty$. Then,  $$\Pi^{Lip}_{p}(X_1,\ldots,X_n; Y)\subset \mathcal{DP}(X_1,\ldots,X_n; Y).$$
\end{proposition}
\begin{pf}
Let $T$ be a Lipschitz $p$-summing multilnear operator. By $(iii)$ in Theorem \ref{thm: equiv p summ multi op}, $T$ factorizes as $T=h_T\circ j_p\circ(i\circ \otimes)$. Consider a sequence $\{(\spaces{x^k}{1}{n})\}_k$ such that $\{(\Tensor{x^k}{1}{n}{})\}_k\subset \Sigma_{\spaces{X}{1}{n}}$ is weakly Cauchy. Since $j_p$ transforms weakly Cauchy sequences into norm convergent sequences,
(see \cite[pp. 40, 49 ]{DJT}) and $h_T$ is norm continuous, then  $\{h_T\circ j_p (\Tensor{x^k}{1}{n}{})\}_k$ is a convergent sequence in $Y$. Thus,
  $\{T(\spaces{x^k}{1}{n}{})\}_k=\{h_T\circ j_p (\Tensor{x^k}{1}{n}{})\}_k$ is a convergent sequence and   we  conclude that  $T$ is Dunford-Pettis.
\end{pf}

 From Theorem \ref{thm: DP} and  Proposition \ref{prop: psum DP} we get:

\begin{corollary}\label{coro: no l1 + p-summ =compact} Let $X_i$ be Banach spaces such that
 $\ell_1\not\subset\CTensor{X}{1}{n}{\pi}$.  Then, for  $1\leq p <\infty$,
    $ {\Pi}_p^{Lip}(X_1,\ldots,X_m; Y)\subset \mathcal{K}(X_1,\ldots,X_m; Y). $
\end{corollary}

\begin{example} Let  $1<q_i<\infty$ such that  $\sum_{i=1}^{n}\frac{1}{q_i}<1$.  Then, for   $1\leq p <\infty$, every
Lipschitz $p$-summing multilinear operator $T\in\mathcal{L}(\spaces{\ell}{q_1}{q_n}; Y)$ is compact. This is because of  Corollary \ref{coro: no l1 + p-summ =compact}  and the fact that   in this case $\CTensor{\ell}{q_1}{q_n}{\pi}$ is reflexive, (see \cite{AriasFarmer96}) .
\end{example}

\section{Hilbert-Schmidt multilinear operators.}\label{sec: Hil-Sch}

 When  the Banach spaces are Hilbert spaces,  absolutely $p$-summing linear operators and Hilbert-Schmidt operators coincide (\cite[Theorem 5.2]{Pelcz66-67}).
   In this section, we determine the relation with Lipschitz $p$-summability.

 In order to remain in a Hilbert-space context, it is necessary to consider the Hilbert-Schmidt tensor product $\otimes_H$, instead of the projective tensor product, since $\CTensor{H}{1}{n}{\pi}$ is  not, in general, a Hilbert space,  while $\CTensor{H}{1}{n}{H}$ is.   We will see briefly how to define  Lipschitz  $p$-sumability when considering other tensor norms.

\subsection{Lipschitz $p$-summability with  respect to other  reasonable crossnorms.}\label{sec: beta Sigma op}
  A  norm $\beta$ on the vector space $\Tensor{X}{1}{n}{}$ is said to be a {\bf reasonable crossnorm} if it has the following  two properties:
   (1) $\beta(x_1\otimes\dots\otimes x_n)\leq \|x_1\|\cdots \|x_n\|$ for every $x_i\in X_i;\; i=1,\ldots n $ and (2)
   For every $x_i^*\in X_i^*$, the linear functional $x_1^*\otimes\dots\otimes x_n^*$ on $\Tensor{X}{1}{n}{}$ is bounded, and $\|x_1^* \otimes\dots\otimes x_n^*\|\leq \|x_1^*\|\cdots \|x_n^*\|$  (see  \cite{DeFlo}, \cite{DieFouSwa}  or   \cite{Ryan-libro}). The main point to pay attention to in  each step is that in this case  we  are working with multilinear functionals whose associated linear operators are continuous on $\CTensor{X}{1}{n}{\beta}$. We will denote the space of such mappings by ${\mathcal{L}_{\beta}\left(\spaces{X}{1}{n}\right)}$.

\begin{definition}\label{def: Lip p-summ beta}
Let  $1 \leq p < \infty$. An $n$-linear operator between Banach spaces    $T\in\mathcal{L}_{\beta}(X_1\ldots , X_n; Y)$ is said to be  Lipschitz $p$-summing with respect to the reasonable crossnorm $\beta$  (briefly, $\beta$-Lipschitz $p$-summing) if there exists $c>0$ such that for $k\in\mathbb{N}$, $i=1,\ldots,k$ and $u_i, v_i \in \cartesian{X}{1}{n}$, the following inequality holds:
 \begin{multline*}
  {\left(\sum_{i=1}^k \left\|{T}\left(u_i\right)-T \left(v_i\right)\right\|^p \right)^{1/p}  }
  \leq   \\ c\cdot \sup \left\lbrace \left(\sum_{i=1}^k \left|{{\varphi}}\left(u_i\right)-{\varphi}\left(v_i\right)\right|^p\right)^{1/p}; \; \varphi \in B_{\mathcal{L}_{\beta}\left(\spaces{X}{1}{n}\right)}\right\rbrace
 \end{multline*}
 \end{definition}

  A Domination and a Factorization Theorems  hold for  Lipschitz $p$-summing multilinear operators with respect to $\beta$  (analogues to (ii) and (iii) in Theorem \ref{thm: equiv p summ multi op}) and many of its properties, such as an inclusion Theorem, can be derived from them.

\subsection{$H$-Lipschitz $p$-summing and   Hilbert-Schmidt multilinear operators coincide.}
In this subsection all the spaces considered $\spaces{H}{1}{n}$, $G$, will be Hilbert spaces and the reasonable crossnorm $\beta$ defined on $\Tensor{H}{1}{n}{}$ will be the Hilbert-Schmidt norm, $\otimes_H$.  In this context, we will say that
$T: \cartesian{H}{1}{n}\rightarrow G$ is $H$-Lipschitz $p$-summing, $T\in \Pi_p^{Lip,H}(\spaces{H}{1}{n};G)$ if
$T$ satisfies Definition \ref{def: Lip p-summ beta} with respect to $\otimes_H$. The best $c$ will be denoted  $\pi_p^{Lip,H}(T)$.

Following \cite[Definition 5.2]{Matos03}, $T$   is said to be Hilbert-Schmidt, $T\in \mathcal{L}_{HS}$,  if
there is an orthonormal basis $(e^k_j)_j$  for each $k = 1, . . . , n$, such that
\begin{align*}\|T\|_{HS}:=\left(\underset{\underset{k=1,..,n}{j_{k}\in J_{k}}}{\sum}\left\Vert T\left(e_{j_{1}}^{1},...,e_{j_{n}}^{n}\right)\right\Vert ^{2}\right)^{\frac{1}{2}} <\infty
\end{align*}
The following  coincidence between classes generalizes the linear result  in \cite{Pelcz66-67}:
\begin{theorem}\label{thm: Hilb-Schm}  For $1\leq p < \infty$, $\mathcal{L_{HS}}(\spaces{H}{1}{n}; G)=\Pi^{Lip,H}_p(\spaces{H}{1}{n}; G)$, where $\left\Vert T\right\Vert _{HS}\leq\pi_{2}^{Lip,H}(T)\leq B_p^n \left\Vert T\right\Vert _{HS}$ and   $B_{p}$ is Khintchine's inequality constant.
\end{theorem}
\begin{pf}
As always,  $\hat{T}$ will denote  the   linear mapping associated to
 ${T}$.
 By  \cite[Proposition 10]{Matos03}, $T\in \mathcal{L_{HS}}(\spaces{H}{1}{n}; G)$ if and only if
the linear mapping $\hat{T}$ is in $ \mathcal{L_{HS}}(\CTensor{H}{1}{n}{H}; G).$  The main result in   \cite{Pelcz66-67} implies  that $\mathcal{L_{HS}}(\CTensor{H}{1}{n}{H}; G)=\Pi_p(\CTensor{H}{1}{n}{H}; G).$ An analogous statement to  Proposition \ref{prop: lineal psumm es Lipsch psumm} holds for the norm $\otimes_H$. Using it  we obtain that if $T\in\mathcal{L}_{HS}$, then
 $T\in \Pi^{Lip,H}_p(\spaces{H}{1}{n}; G)$.

 To prove  the reverse  inclusion $ \Pi^{Lip,H}_p\subset \mathcal{L}_{HS}$, we consider first  the case $p=2$.
Let $T\in\Pi^{Lip,H}_{2}\left(H_{1},H_{2},\ldots,H_{m};G\right) $  and  let
 $\left(e_{i_{k}}^{k}\right)_{i_{k}\in I_{k}}$ be an orthonormal
basis for $H_{k},$ $k=1,...,n$. Fix $J_{k}\subseteq I_{k},$ an
arbitrary finite subset of $I_{k}.$ From  Definition \ref{def: Lip p-summ beta} we have that in this particular case
\begin{align*}
\left(\underset{\underset{k=1,..,n}{j_{k}\in J_{k}}}{\sum}\left\Vert T\left(e_{j_{1}}^{1},...,e_{j_{n}}^{n}\right)\right\Vert ^{2}\right)^{\frac{1}{2}} & \leq\pi_{2}^{Lip,H}\left(T\right)\underset{\varphi\in B_{\mathcal{L}_{HS}}}{\sup}\left(\underset{\underset{k=1,..,n}{j_{k}\in J_{k}}}{\sum}\left|\varphi\left(e_{j_{1}}^{1},...,e_{j_{n}}^{n}\right)\right|^{2}\right)^{\frac{1}{2}}\\
 & =\pi_{2}^{Lip,H}\left(T\right).
\end{align*}
Therefore $\left\Vert T\right\Vert _{HS}\leq\pi_{2}^{Lip,H}\left(T\right).$

For  $1\leq p \leq 2$ we get, from this case and from   Proposition \ref{prop: multi usual p,q contentions} stated   with the norm $\otimes_H$  (namely, $\Pi^{Lip,H}_p\subset \Pi^{Lip,H}_2$),  that  $ \Pi^{Lip,H}_p\subset \mathcal{L}_{HS}$ and $\left\Vert T\right\Vert _{HS}\leq\pi_{2}^{Lip,H}\left(T\right).$

 Finally, let  $p>2$. We will see that whenever $T:H_{1}\times H_{2}\times\cdots\times H_{n}\rightarrow G$ is an
element of $\Pi^{Lip,H}_{p}\left(H_{1},H_{2},\ldots,H_{n};G\right)$, then  $\| T\|_{HS}  \leq\left(B_{p}\right)^{n}\pi_{2}^{Lip,H}\left(T\right)$.

   Let  $\left(e_{i_{k}}^{i}\right)_{i_{k}\in I_{i}}\subseteq H_{i}$ be
an orthonormal basis of $H_{i},$ $i=1,...,n$. For every finite
subset $J_{i}\subseteq I_{i},$ with $m_{i}$ elements, we consider
$\left(e_{j_{k}^i}^{i}\right)_{j_{k}^i\in J_{i}}$ ordered linearly and for each $i=1,..,n,$ write $e_{j^i_1}^{i},...,e_{j^i_{m_{i}}}^{i}.$   If $\{r_n\}_{n}$ are  the Rademacher functions and writing   $dt$ for $dt_1\cdots dt_n$ and $d\mu$ for $d\mu(\varphi)$, we have
\begin{align*}
  & \left(\overset{m_{l}}{\underset{\underset{l=1,..,n}{{k_l}=1}}{\sum}}\left\Vert T\left(e_{j^1_{k_1}}^{1},...,e_{j^n_{k_n}}^{n}\right)\right\Vert ^{2}\right)^{\frac{1}{2}} =\\
 & \left(\underset{\left[0,1\right]^{n}}{\idotsint}\left\Vert T\left(\overset{m_{1}}{\underset{k_{1}=1}{\sum}}r_{k_{1}}\left(t_{1}\right)e_{j_{k_{1}}}^{1},...,\overset{m_{n}}{\underset{k_{n}=1}{\sum}}r_{k_{n}}\left(t_{n}\right)e_{j_{k_{n}}}^{n}\right)\right\Vert ^{2}dt_{1}\cdots dt_{n}\right)^{\frac{1}{2}}\leq \\
 &  \left(\underset{\left[0,1\right]^{n}}{\idotsint}\left\Vert T\left(\overset{m_{1}}{\underset{k_{1}=1}{\sum}}r_{k_{1}}\left(t_{1}\right)e_{j^1_{k_1}}^{1},...,\overset{m_{n}}{\underset{k_{n}=1}{\sum}}r_{k_{n}}\left(t_{n}\right)e_{j^n_{k_n}}^{n}\right)\right\Vert ^{p}dt_{1}\cdots dt_{n}\right)^{\frac{1}{p}}\leq \\
&\pi^{Lip,H}_{p} \left(T\right) \left(\underset{\left[0,1\right]^{n}}{\idotsint}
\underset{B_{\mathcal{L}_{HS}}}{\int}\left|\varphi\left(\overset{m_{1}}{\underset{k_{1}=1}
{\sum}}r_{k_{1}}\left(t_{1}\right)e_{j^1_{k_1}}^{1},...,\overset{m_{n}}{\underset{k_{n}=1}
{\sum}}r_{k_{n}}\left(t_{n}\right)e_{j^n_{k_n}}^{n}\right)\right|^{p}d\mu dt\right)^{\frac{1}{p}} \leq\\
 &\left(B_{p}\right)^{n}\pi^{Lip,H}_p\left(T\right)\left(\underset{B_{\mathcal{L}_{HS}}}{\int}\left(\overset{m_{l}}{\underset{\underset{l=1,..,n}{k_{l}=1}}{\sum}}\left|\varphi\left(e_{j^1_{k_1}}^{1},...,e_{j^n_{k_n}}^{n}\right)\right|^{2}\right)^{\frac{p}{2}}d\mu\left(\varphi\right)\right)^{\frac{1}{p}}  \leq \\ & \left(B_{p}\right)^{n}\pi^{Lip,H}_p\left(T\right).
\end{align*}
\end{pf}

\section[A Chevet-Saphar type norm]{Tensorial approach to Lipschitz  $p$-summing multilinear operators.}\label{sec: tensors}

Here we study the tensorial side of the $p$-summability of multilinear mappings. By one hand, we will construct a tensor normed space  in such a way that its dual
 is  isometric to  the space of   Lipschitz $p$-summing multilinear operators.   On the other hand, Lipschitz $p$-summing multilinear operators will be characterized  in terms of a Lipschitz  condition of a mapping of the form $Id\otimes T$.

 To prove the first statement, we  consider  decompositions of the elements  $z\in \Tensor{X}{1}{n}{}\otimes Y$ of the following form:
\begin{equation}\label{eq: Angulo decomposition}
 z=\sum_{i=1}^{k}(p_i-q_i)\otimes y_i \;\;\; {\mbox{where }}\; p_i, q_i\in \Sigma_{\spaces{X}{1}{n}},\; y_i\in Y.
\end{equation}

Given a pair of conjugate indexes  $1\leq p, p'\leq \infty$ and a  $z\in \Tensor{X}{1}{n}{}\otimes Y$ we define
$$d^{Lip}_p(z):=\inf \left\lbrace \left( \sup_{\varphi\in B_{{\mathcal{L}}({\Sigma_{X_1,\ldots, X_n}})}}\left\lbrace(\sum_{i=1}^{n}|\varphi(p_i)-\varphi(q_i)|^{p'})^{\frac{1}{{p'}}} \right\rbrace\right)\cdot\|(y_i)_{i=1}^n\|_p\right\rbrace$$
taking the infimum over all representations of the form (\ref{eq: Angulo decomposition}).

It can be proved that  $d^{Lip}_p$ defines a reasonable crossnorm on
 $\Tensor{X}{1}{n}{}\otimes Y$.  {This definition generalizes  the (right) Chevet-Saphar tensor norm $d_p$ defined on spaces of the form $X\otimes Y$ (see \cite[(6.5)]{Ryan-libro}).} The special form of the representations
  (\ref{eq: Angulo decomposition})  of an arbitrary tensor will
   capture the Lipschitz character of the  $\Sigma$-operator associated to
   Lipschitz $p$-summing multilinear operators.
  {The completion of the tensor space with respect to the norm
   $d_p^{Lip}$ will be denoted $(X_1\hat{\otimes}\ldots\hat{\otimes} X_n\hat{\otimes} Y, d_p^{Lip})$. }
\begin{theorem}\label{dual representation} Let  $p, p'$ be  conjugate indexes,  $1<p \leq \infty$ and  $\spaces{X}{1}{n}, Y$  Banach spaces.  Then
$$(\CTensor{X}{1}{n}{}\hat{\otimes}Y,{d^{Lip}_p})^*=\Pi^{Lip}_{p'}(\spaces{X}{1}{n}; Y^*). $$
\end{theorem}
  This  result can be proved by  generalizing in a natural way the linear proof.

In the linear setting,  the absolutely $p$-summability of a linear operator
$T$ is equivalent to the   boundedness  of the operator
$1\otimes T:\ell_p\otimes_{\varepsilon} X\rightarrow \ell_p\otimes_{\Delta_p} Y$  { where $\epsilon$ denotes the injective tensor norm  and $\Delta_p$ is the norm on $\ell_p\otimes Y $ induced  by the embedding $\ell_p\otimes Y \hookrightarrow \ell_p[Y]$ (see \cite[4.1, 7.1]{DeFlo})}. For multilinear operators, such  characterization reads as follows:

\begin{theorem}\label{thm: equiv p-Lips norm} Let $\spaces{X}{1}{n}, Y$ be Banach spaces and $ 1 \leq p  < \infty$. The following conditions are equivalent for a bounded multilinear mapping $T$:
\begin{itemize}
  \item[(i)]  $T\in\mathcal{L}(\spaces{X}{1}{n}; Y)$ is Lipschitz $p$-summing.
  \item[(ii)] There exists  $C>0$ such that
  $$\Delta_p\left((1\otimes \hat{T})\left(\sum_{i=1}^{k} e_i\otimes (a_i -b_i  )\right)\right)\leq C \cdot \varepsilon\left(\sum_{i=1}^{k}  e_i\otimes (a_i-b_i)\right)  $$
   {where}  $a_i, b_i \in \Sigma_{\spaces{X}{1}{n}}$, $(e_i)_i$ denotes the canonical basis of $\ell_p$ and $$1\otimes \hat{T}:\ell_p\otimes
   (\CTensor{X}{1}{n}{\pi})\rightarrow \ell_p\otimes Y. $$
  \end{itemize}
In this case ${\pi}_p^{Lip}(T)=\inf C$ in $\mathit{(ii)}$.
  \end{theorem}
\begin{pf}
$(i) \Rightarrow (ii)$:  Let $f_T$ be the $\Sigma$-operator associated to $T$. Writing $\mathcal{L}:=\mathcal{L}(\Sigma_{X_1,\ldots, X_n})$, we have
\begin{multline*}
\Delta_p\left(\sum e_i\otimes (f_T(a_i) -f_T(b_i)  )\right)  =  \left( \sum \| f_T(a_i)-f_T(b_i) \|^p\right)^{1/p}\leq \\
     \pi_p^{Lip}(T)  \sup\limits_{\varphi\in B_{\mathcal{L}}} \left(\sum |\varphi(a_i)-\varphi(b_i)|^p\right)^{1/p}
   = \\ \pi_p^{Lip}(T)\sup_{\varphi\in B_{\mathcal{L}}}\sup_{\mu\in B_{\ell_{p'}}} \left|\sum \mu_i (\varphi(a_i)-\varphi(b_i)) \right|
   =  \pi_p^{Lip}(T)\varepsilon\left(\sum e_i\otimes (a_i-b_i)\right).
\end{multline*}

$(ii) \Rightarrow (i)$: Let  $(a_i), (b_i) \in \Sigma_{\spaces{X}{1}{n}}$.   Then,
\begin{multline*}
\left(\sum  \|f_T(a_i)-f_T(b_i) \|^p\right)^{1/p}=\Delta_p\left(\sum e_i\otimes (f_T(a_i) -f_T(b_i)  )\right)\\
  \leq  C \cdot \varepsilon\left(\sum e_i\otimes (a_i-b_i  )\right)
    =C \sup\limits_{\varphi\in B_{\mathcal{L}},\mu\in B_{\ell_{p'}}}  \left|\sum \mu_i (\varphi(a_i)-\varphi(b_i)  )\right|\\
  =  C\sup\limits_{\varphi\in B_{\mathcal{L}}} \left(\sum |\varphi(a_i)-\varphi(b_i)|^p\right)^{1/p}.
\end{multline*}
\end{pf}

Notice that  $Id \otimes \hat{T}$ is bounded (or Lipschitz) on all of its domain exactly when $\hat{T}$  is a linear $p$-summing operator.

\section{Other non-linear  notions of $p$-summability}\label{sec: comparison}

\subsection{Relation with other multilinear generalizations of $p$-summability}\label{subs: comparison}

    In Proposition \ref{prop: lineal psumm es Lipsch psumm}, we have already exposed the relation with two  generalizations of absolutely $p$-summing operators that have  appeared in the literature. Another  class is that of the so called fully $p$-summing   \cite{Matos03}  or multiple $p$-summing \cite{BomPerVil} multilinear operators  (both classes coincide).  There is not a general containment relation with Lipschitz $p$-summing multilinear operators. As in the cases  of the absolutely $p$-summing multilinear operators and the $np$-dominated operators introduced by A. Pietsch in \cite{Pietsch},  it happens that these  classes do not admit a characterization in terms of a   Pietsch-type factorization,  since  there are scalar valued multilinear operators that do not belong to such classes (see \cite{PV} and \cite{Pietsch}).  This fact makes these notions essentially different  to  the one introduced here (see Proposition \ref{prop: multi-ideal} and Theorem \ref{thm: equiv p summ multi op}).

\subsection{{Relation with the notion of $p$-summability for Lipschitz mappings between metric spaces}}\label{subs: Lips}

 In \cite{FJ09}, J.D. Farmer and W.B. Johnson introduced  the notion of Lipschitz $p$-summability  for  mappings between metric spaces, in  the following way: ``The Lipschitz $p$-summing  ($1\leq p <\infty$) norm, $\pi_p^L(T)$, of a (possibly nonlinear) mapping $T:X\rightarrow Y$  between metric spaces   is the smallest constant $C$ so that for all $(x_i)_i$, $(y_i)_i$ in $X$ and all positive reals $a_i$,
 $$\sum a_i\|Tx_i-Ty_i\|^p\leq C^p \sup_{f\in B_{X^{\#}}} \sum a_i\|fx_i-fy_i\|^p.$$
 Here $B_{X^{\#}}$ is the unit ball of ${X^{\#}}$, the Lipschitz dual of $X$''. As they noted,   the definition is the same if one considers $a_i=1$.
 In  \cite[Theorem 2]{FJ09} the authors  prove that  { this notion coincides with the notion of absolutely $p$-summing operators when the mappings are linear operators.}  We have already mentioned that $\Sigma$-operators are   Lipschitz mappings between metric spaces. Then, it is natural to ask  if the $\Sigma$-operators   that are  $p$-summing  mappings in the sense  of \cite{FJ09}  (i.e. those that satisfy the inequality above),  are exactly  the  absolutely $p$-summing $\Sigma$-operators as defined in  (\ref{eq: def sigma p sumante}).   Directly  from the definitions we see that every absolutely   $p$-summing $\Sigma$-operator    is $p$-summing in the sense of \cite{FJ09}. So,
  \begin{question}    Do both notions coincide for  $\Sigma$-operators?
 \end{question}

   \subsection{Lipschitz $p$-summing  polynomials}
  A mapping $P:X\rightarrow Y$ between Banach spaces is a {\sl homogeneous poylnomial of degree $n$}
if there exists a multilinear mapping $T_P: \cartesian{X}{}{}\rightarrow Y$ such that $P(x)=T_P(x,\stackrel{n}{\ldots},x)$ (see \cite{Dineen book}).
Let  $\Poll{n}{X}{Y}$ (resp. $\Pol{n}{X}$) denote the
 space of $n$-homogeneous bounded  polynomials from $X$ to $Y$ (resp. from $X$ to the scalar field), normed with the supremum norm.
 The mappings $j_p,i_p$ and $i_Y$ will denote the analogues  to those described before Theorem \ref{thm: equiv p summ multi op}. The mapping $i_X^n$  is the $n$-homogeneous polynomial  defined as $i_X^n(x)=x\otimes\stackrel{n}{\cdots}\otimes x$. Adapting the arguments already used to prove Theorem \ref{thm: equiv p summ multi op}, it is possible to prove the following result.

\begin{theorem}\label{thm: equiv p summ pol} Let $1 \leq p < \infty$ and  let $P\in\Poll{n}{X}{Y}$ be an  $n$-homogeneous polynomial  between Banach spaces. The following conditions for $P$ are equivalent:
\begin{enumerate}
\item   There exists $c>0$ such that for $k\in\mathbb{N}$, $i=1,\ldots,k$ and $u_i, v_i \in X$,
  \begin{multline*} \label{defsumantes.1}
  {\left(\sum_{i=1}^k \left\|{P}\left(u_i\right)-P \left(v_i\right)\right\|^p \right)^{1/p}  }
  \leq  \\  c\cdot \sup \left\lbrace \left(\sum_{i=1}^k \left|{{\varphi}}\left(u_i\right)-{\varphi}\left(v_i\right)\right|^p\right)^{1/p}; \; \varphi \in B_{\Pol{n}{X}}\right\rbrace.
 \end{multline*}
 \item[]

\item There is a constant $c\geq0$ and a regular probability measure $\mu$ on $\left( B_{\Pol{n}{X}},w^*\right)$ such that for each $u, v \in X$ we have that
\begin{equation*} \label{eqn: dominacion}
{\left\|P\left(u\right)-P\left(v\right)\right\| }
 \leq  c\cdot\left(\int_{B_{\Pol{n}{X}}} \left|\varphi\left(u\right)-\varphi\left(v\right)\right|^p\,d\mu(\varphi)\right)_{\ .}^{1/p}
\end{equation*}
\item  There exist a regular Borel probability measure $\mu$ on the space $\left(B_{\Pol{n}{X}},w^*\right)$, a  subset $X_p:=(j_p\circ i_X^n) (X) \subset L_p\left(\mu\right)$ and a Lipschitz function ${{h_P}}: {X_p}\rightarrow {Y}$ such that
 $ h_P\circ j_p \circ i_X^n = P$.

     \[
\begin{tikzcd}\label{p-summ diag4}
        X \arrow{r}{P}\arrow{d}[swap]{i_X^n} &  Y\\
        i_X^n(X)  \arrow{r}{{j_p}_{|_{}}}\arrow[phantom,"\cap"]{d} &  X_p \arrow[phantom,"\cap"]{d}\arrow{u}[swap]{h_P}\\
        C(B_{\Pol{n}{X}},w^*)\arrow{r}{{j_p}}   &    L_p\left(\mu\right)\\
 \end{tikzcd}
\]

\item There exist  a  probability space $(\Omega,\Sigma,\mu)$, a $n$-homogeneous polynomial
$\nu: X \rightarrow L_{\infty}(\mu)$ and   a Lipschitz
function $\funcion{\tilde{h_P}}{L_p\left(\mu\right)}{\ell_{\infty}^{B_{Y^*}}}$
such that $i_Y\circ P=\hat{h}_P\circ i_p\circ \nu$.
the following diagram commutes:
     \[
\begin{tikzcd}\label{p-summ poly linfty}
       X \arrow{r}{P}\arrow{dd}[swap]{\nu} &  Y \arrow{rd}{i_Y} & \\
        & & \ell_{\infty}^{B_{Y^*}} \\
        L_{\infty}(\mu)\arrow{r}{{i_p}}   &    L_p\left(\mu\right)\arrow{ru}{\tilde{h_P}} & \\
 \end{tikzcd}
\]
\end{enumerate}
If $\pi_p^{Lip}(P):=\inf\{c;\; \mathit{(1)} \; \mbox{holds}\}$,
   then  $\pi_p^{Lip}(P)  =\inf\{c;\; \mathit{(2)} \; \mbox{holds}\}$ and   $\normasub{\tilde{h_P}}{Lip}= \pi_p^{Lip}(P)=\normasub{{h_P}}{Lip}$ when the spaces are  real  and  $\pi_p^{Lip}({P}) = \normasub{{h_P}}{Lip} \leq \normasub{\tilde{h_P}}{Lip}\leq \sqrt{2}\pi_p^{Lip}(P)$ when the spaces are complex.
\end{theorem}

 \begin{definition}\label{def p summing poly} An $n$-homogeneous polynomial  $P \in\mathcal{P}(^n X; Y)$ is {\bf Lipschitz $p$-summing} if and  only it satisfies any of the equivalent conditions in Theorem \ref{thm: equiv p summ pol}.
\end{definition}
Important  properties of these polynomials can be
derived from  Theorem \ref{thm: equiv p summ pol}. This is the case of the Inclusion Theorem and the fact that   whenever the symmetric projective tensor product contains no copies of $\ell_1$, then every  Lipschitz $p$-summing polynomial is compact.

The  following   homogenous polynomials, which  are analogues to some of the  multilinear examples given in Section \ref{sec: Lips multi}, are Lipschitz $p$-summing.
This fact can be  proved in like manner as the proof of the linear case $n=1$ \cite[Examples 2.9]{DJT}.

 Let $K$  be a compact Hausdorff space and  $\mu$ a positive regular Borel measure      on $K$. If   $h\in L_p(\mu)$,  the  $n$-homogeneous polynomial  $P_h\in\mathcal{P}({C(K)}; L_p(\mu))$ defined as $P_h(g):=h(w)\cdot (g(w))^n$
is  Lipschitz $p$-summing.

Given a sequence  $\lambda=(\lambda_k)_{k}\in \ell_p$, the
 diagonal polynomial
  $P_{\lambda}\in\mathcal{P}({\ell_{\infty}};  \ell_p)$, defined as
$   P_{\lambda}\left((a_k)_{k}\right):=({\lambda_k\cdot a_k^n})_{k}$, is Lipschitz $p$-summing.

\section*{Acknowledgment}

\noindent
The authors want to thank  Luisa F. Higueras-Monta\~no and Samuel Garc\'\i a-Hern\'andez for their contributions especially  in Sections \ref{sec: Hil-Sch} and \ref{sec: tensors}, respectively.


\begin{thebibliography}{99}



 \bibitem{AlbBayPelSeo}Albuquerque, N.; Bayart, F.; Pellegrino, D.; Seoane-Sep\'ulveda, J. B. Optimal Hardy-Littlewood type inequalities for polynomials and multilinear operators, Israel J. Math. 211 (2016), no. 1, 197-220.

\bibitem{AriasFarmer96} Arias, Alvaro; Farmer, Jeff D. On the structure of tensor products of $l_p$-spaces, Pacific J. Math. 175 (1996), no. 1, 13-37.




  \bibitem{BenyLind}  Benyamini Y.;  Lindenstrauss J. {Geometric Nonlinear Functional Analysis Volume 1}; American Mathematical Society Colloquium Publications, 48, 2000.

 \bibitem{BlBoPelRu} Blasco, O.; Botelho, G.; Pellegrino, D.; Rueda, P. Summability of multilinear mappings: Littlewood, Orlicz and beyond, Monatsh. Math. 163 (2011), no. 2, 131-147.

 \bibitem{Bombal90} Bombal, F.; On some properties of Banach spaces.  Rev. Real Acad. Cienc.
Exact. F\' \i s. Natur. Madrid, 84 (1990), no. 1, 83-116.


  \bibitem{BomPerVil} Bombal F.; P\'erez-Garc\'ia, D.;  Villanueva I. {Multilinear extensions of Grothendieck's Theorem}, Q. J. Math. \textbf{55} (2004), no. 4,  441-450.


 \bibitem{BotPelRue2}  Botelho G.; Pellegrino D.;  Rueda P. { A unified Pietsch domination theorem, J. Math. Anal. Appl}. {365} (2010), no. 1, 269-276.

\bibitem{CaliPelle} \c{C}ali\c{s}kan, E.; Pellegrino, D.  On the multilinear generalizations of the concept of absolutely summing operators,  Rocky Mountain J. Math. 37 (2007), no. 4, 1137-1154.

\bibitem{CarDimant03}  Carando D.;  Dimant V.  On summability of bilinear operators, Math. Nachr. 259 (2003), 3-11.




\bibitem{DeFlo} Defant, A.; Floret, K. Tensor norms and operator ideals. \textit{North-Holland Mathematics Studies}, 176, Amsterdam, 1993.

\bibitem{DefSev} Defant A.;   Sevilla-Peris P.  A new multilinear insight on Littlewood's $\frac{4}{3}$-inequality. J. Funct. Anal. 256 (2009), no. 5, 1642-1664.

 \bibitem{DJT}  Diestel J.;  Jarchow, H.;  Tonge A. {Absolutely summing operators. Cambridge Univ. Press}, 1995.

\bibitem{DieFouSwa} Diestel, J.; Fourie, J. H.; Swart, J.  The metric theory of tensor products. Grothendieck's r\'esum\'e revisited. \textit{American Mathematical Society, Providence, RI}, 2008.


\bibitem{Dimant03}  Dimant V. {Strongly $p$-summing multilinear operators}; J. Math. Anal. Appl. {278} (2003), 182-193.

\bibitem{DimSev}  Dimant, V.; Sevilla-Peris, P. Summation of coefficients of polynomials on $\ell_p$  spaces. Publ. Mat. 60 (2016), no. 2, 289-310.

  \bibitem{Dineen book} Dineen, S.  Complex analysis on infinite-dimensional spaces. Springer Monographs in Mathematics. Springer-Verlag London, Ltd., London, 1999.







\bibitem{FJ09}  Farmer, J.D.; Johnson, W. B.
Lipschitz $p$-summing operators.
Proc. Amer. Math. Soc. 137 (2009), no. 9, 2989-2995.


\bibitem{MFU2001} Fern\'andez-Unzueta, M.  Dunford-Pettis and Dieudonn\'e polynomials on Banach spaces, Illinois J. Math. 45 (2001), no. 1, 291-307.

\bibitem{MFU} Fern\'andez-Unzueta, M. The Segre cone of Banach spaces and multilinear mappings, Linear and Multilinear Algebra, 68 (2020), no. 3, 575-593.

\bibitem{F-UG-H}  Fern\'andez-Unzueta, M.;  Garc\'\i a-Hern\'andez, S.  {Multilinear operators  factoring through Hilbert spaces}. Banach J.  Math. Anal. 13 (2019), no. 2, 234-254.


\bibitem{F-UG-H 2}  Fern\'andez-Unzueta, M.;   Garc\'\i a-Hern\'andez, S.  $(p,q)$-dominated multilinear operators and Laprest\'e tensor norms. J.  Math. Anal.  App. 470 (2019), no. 2, 982-1003.



 \bibitem{LP} Lindenstrauss J.;  Pe\l czy\'nski A. {Absolutely summing operators in $\mathcal{L}_p$ and their applications}. Studia Math. 29 (1968), 275-326.

   \bibitem{Matos03} Matos M.C. {Fully absolutely summing  and Hilbert-Schmidt multilinear mappings}, Collect. Math.  54 (2003), no.2 111-136.


 \bibitem{Pelcz66-67} Pe\l czy\'nski, A. A characterization of Hilbert-Schmidt operators. Studia Math. 28 (1966/67) 355-360.


   \bibitem{PelRib}  Pellegrino D.;  Ribeiro, J. {On multi-ideals and polynomial ideals of Banach spaces: a new approach to coherence and compatibility}, Monatsh Math.  Monatsh. Math. 173 (2014), no. 3, 379-415.



  \bibitem{Per-Gar}   P\'erez-Garc\'\i a D. {Comparing different classes of absolutely summing multilinear operators}; Arch. Math. {85} (2005), 258-267.

  \bibitem{PV}  P\'erez-Garc\'\i a D.; Villanueva I. {Multiple summing operators on Banach Spaces}, J. Math. Anal. Appl. {285} (2003) no. 1, 86-96.




  \bibitem{Pietsch} Pietsch A. {Ideals of multilinear functionals (designs of a theory)}, Proceedings of the second international conference on operator algebras, ideals, and their applications in theorical physics (Leipzig, 1983), Teubner-Texte Math., vol. 67, Teubner, Leipzig, 1984, 185-199.


  \bibitem{Pisier}   Pisier G. {Factorization of linear operators and geometry of Banach spaces}. CBMS Regional Conference Series in Mathematics, 60.
      American Mathematical Society, Providence, RI, 1986.



\bibitem{Pisier Non-comm p-sum} Pisier G.  Non-commutative vector valued Lp-spaces and completely p-summing maps.
Ast\'erisque No. 247 (1998).

\bibitem{Ryan-libro}  Ryan, R.  Introduction to tensor products of Banach spaces. Springer Monographs in Mathematics. \textit{Springer-Verlag London, Ltd., London}, 2002.




 \end{thebibliography}
\end{document}